%% file: dirac-massive20130406.tex
\documentstyle[12pt,bezier]{article}

\title{Global solution to a cubic nonlinear  Dirac  equation  in $1+1$ dimensions}
\author{Yongqian Zhang
\\ {\small School of Mathematical Sciences, Fudan University,
Shanghai 200433} \\ {\small Key Laboratory of Mathematics for
Nonlinear Sciences }
\\ {\small Email address: yongqianz@fudan.edu.cn} }

\begin{document}
\maketitle

\makeatletter
\renewcommand{\theequation}{\thesection.\arabic{equation}}
\@addtoreset{equation}{section}
\linespread{2}\selectfont
\makeatother

\newtheorem{lemma}{Lemma}[section]
\newtheorem{proposition}{Proposition}[section]
\newtheorem{theorem}{Theorem}[section]
\newtheorem{definition}{Definition}[section]
\newtheorem{remark}{Remark}[section]

\begin{abstract}This paper studies a class of nonlinear  Dirac equations with cubic terms in $R^{1+1}$, which  include the equations for the massive Thirring model and the massive Gross-Neveu model. Under the assumptions that the initial data has small charge, the global existence  of the solution in $H^1$ are proved. The proof is given by introducing some Bony functional to get the uniform estimates on the nonlinear terms and the uniform bounds on the local smooth solution, which enable us to extend the local solution globally in time. Then $L^2$-stability estimates for these solutions are also established by a Lyapunov functional and the global existence of weak solution in $L^2$ is  obtained.
\end{abstract}

\section{Introduction}

We consider the nonlinear  Dirac equations
\begin{equation}\label{eq-dirac}
\left\{ \begin{array}{l} i(u_t+u_x)=-mv +N_1(u,v), \\i( v_t-v_x)=-mu+N_2(u,v),
\end{array}
\right.
\end{equation}
with initial data
\begin{equation}\label{eq-dirac-initialv}
(u, v)|_{t=0}=(u_0(x), v_0(x))
\end{equation}
where $(t,x)\in R^2$, $(u,v)\in \mathbf{C}^2$, $m\ge0$. The nonlinear terms   satisfy the following:
\begin{description}
\item[(A1)] $N_1$ and $N_2$ have following form:\[N_1(u,v)=(\alpha_1 u+ \alpha_2 \overline{u})(\alpha_3 |v|^2+\alpha_4 v^2 +\alpha_5 \overline{v}^2),\]
    \[N_2(u,v)=(\beta_1 v+ \beta_2 \overline{v})(\beta_3 |u|^2+\beta_4 u^2 +\beta_5 \overline{u}^2),\]
    where $\alpha_k$ and $\beta_k$ $(k=1,\cdots, 5)$ are constants.
\item[(A2)] For any $(u,v)$, \[ \Re \big(i \overline{u}N_1(u,v) +i \overline{v}N_2(u,v)\big)=0.\]Here and in sequel $\Re Z$ denotes the real part of $Z \in \mathbf{C}$.
\end{description}
Many physical models like massive Thirrring model and massive Gross-Neveu model verify (A1) and (A2). Indeed,
for  massive Thirrring model (\cite{thirring}), we have
\[ (N_1,N_2)=\nabla_{(\overline{u},\overline{v})} \alpha |u|^2|v|^2=(\alpha u|v|^2, \alpha v|u|^2), \, \alpha\in R^1,\] therefore
\[ \Re \big(i \overline{u}N_1 \big)=\Re \big(i \overline{v}N_2 \big)=0.\]
For the massive Gross-Neveu model (\cite{gross-neveu}), we have
\[ (N_1, N_2)=\nabla_{(\overline{u}, \overline{v})} \alpha\big(\overline{u}v +u\overline{v} \big)^2=\big(2\alpha v(\overline{u}v +u\overline{v}),2\alpha u (\overline{u}v +u\overline{v})\big), \, \alpha\in R^1,\]
then
\[ \Re \big(i \overline{u}N_1 +i \overline{v}N_2\big)=\Re\{2i\alpha\big(\overline{u}v +u\overline{v} \big)^2\}= 0.\]

The nonlinear Dirac equation is important in quantum mechanics and general relativity (\cite{gross-neveu} and \cite{thirring}). There are a number of works on the local and global well-posedness of Cauchy problem in different Sobolev spaces, see for examples \cite{bachelot}, \cite{bachelot-2}, \cite{bournaveas-zouraris}, \cite{cacciafesta1}, \cite{cacciafesta2}, \cite{candy},  \cite{deldado},  \cite{dias-figueira1}, \cite{dias-figueira2}, \cite{escobedo}, \cite{Esteban-Lewin-Sere}, \cite{huh}, \cite{machihara}, \cite{pelinovsky}, \cite{selberg}, \cite{zhang} and the references therein. The survey of the well-posedness in the nonlinear Dirac equation in one dimension is given in \cite{pelinovsky}.
In this paper we consider a class of cubic nonlinear Dirac equations in one dimension, which include the equation for Massive Gross-Neveu model as an example. To our knowledge, the global existence of  solution in $H^1$ or in $L^2$  for the Dirac equation of Massive Gross-Neveu model is still open \cite{pelinovsky}.
And it was pointed out by Pelinovsky \cite{pelinovsky} that the apriori estimates of the $L^p-$norm $(||u(t)||_{L^p}^p+||v(t)||_{L^p}^p)^{1/p}$ for the equations like Gross-Neveu model include  nonlinear terms, which may lead to the blow-up of the $L^{\infty}$ and $H^1$ norms. In addition, the energy of Dirac equation is non-positivity . We will overcome the difficulties caused by  the nonlinear terms in the apriori estimates for $p=2$ and the non-positivity  of the energy by introducing a Bony type functional $Q(t)=\int\int_{x<y} |u(t,x)|^2|v(t,y)|^2dxdy$, which is similar to the Glimm's interaction potential for conservation laws (\cite{glimm}) and is also similar to the Bony functional for Boltzmann  equations (\cite{bony} and \cite{ha-tzavaras}). Then, following Bony's approach \cite{bony},  for any local $H^1-$solution $(u,v)$ with small initial charge $\int_{-\infty}^{\infty} (|u_0|^2+|v_0|^2)dx$, we can get the control on $\int^t_0 \int^{\infty}_{-\infty} |u(x,\tau)|^2|v(x,\tau)|^2 dx d\tau$ by using the conservation of the charge $\int_{-\infty}^{\infty} (|u(t,x)|^2+|v(t,x)|^2)dx$; see Lemma \ref{lemma-bony-functional}. Here the conservation of charge comes from the assumption (A2) on the nonlinear terms. These bounds enable us to apply the characteristic method  to the equations (\ref{eq-dirac1}) to get the uniform $L^{\infty}$ bounds on $|u|^2$ and $|v|^2$  for  extending the local $H^1-$solution globally in time; see Lemmas \ref{lemma-integral} and \ref{lemma-bound} for the $L^{\infty}$ bounds.
 Then, to study the stability of the $H^1-$solutions, we have (\ref{eq-dirac-difference3}) and introduce Lyapunov functional $L_1(t)+KQ_1(t)$ to derive the $L^2$ estimates on the difference between these $H^1-$solutions. Finally global weak solution in $L^2$ is obtained.
 We remark that Glimm interaction potential was first used by Glimm \cite{glimm} to establish the global existence of the BV solution for the system of conservation laws, then was used to study the general conservation laws, and that Bony type functional was used to study Boltzmann equations and wave maps in $R^{1+1}$, see for instance \cite{glimm},  \cite{bressan}, \cite{dafermos}, \cite{tartar}, \cite{bony},   \cite{ha-tzavaras}and \cite{zhou} and the references therein.

The main results is stated as follows.
\begin{theorem}\label{thm-1-existence}
Suppose that $(u_0,v_0)\in H^1(R^1)$ and that the norm $||u_0||_{L^2(R^1)}^2 +||v_0||_{L^2(R^1)}^2\le C_0$  for some small constant $C_0>0$. Then (\ref{eq-dirac}-\ref{eq-dirac-initialv}) admits a unique global solution $(u,v)\in C([0,\infty); H^1(R^1))\cap C^1([0,\infty); L^2(R^1))$.
\end{theorem}

\begin{theorem}\label{thm-2-convergence} Suppose that $(u_{0k},v_{0k})\in H^1(R^1)$ with $||u_{0k}||_{L^2(R^1)}^2 +||v_{0k}||_{L^2(R^1)}^2\le C_0$ $ (1\le k< \infty)$ such that $||u_{0k}-u_{0\infty}||_{L^2(R^1)} +||v_{0k}-v_{0\infty}||_{L^2(R^1)}\to 0$ as $k\to \infty$ for some $(u_{0\infty},v_{0\infty})\in L^2(R^1)$. Let $(u_k,v_k)$ be the classical solution of (\ref{eq-dirac}) $(1\le k<\infty)$ taking $(u_{0k},v_{0k})$ as its initial data. Then there exists a weak solution $(u_{\infty},v_{\infty})\in C([0,+\infty),L^2(R^2))$ of (\ref{eq-dirac}) such that
\[ ||u_k-u_{\infty}||_{C([0,T],L^2(R^2))}+||v_k-v_{\infty}||_{C([0,T],L^2(R^2))}\to 0\]
as $k\to \infty$ for any $T\ge 0$.
\end{theorem}

\section{Global classical solution }

 Since $C^{\infty}_c(R^1)$ is dense in $H^1(R^1)$, we consider in this section the problem for $(u_0,v_0)\in C^{\infty}_c(R^1)$. (\ref{eq-dirac}) is a semilinear hyperbolic system. Local existence of the solutions in $H^1$ for the semilinear strictly hyperbolic system like (\ref{eq-dirac}) can be proved by the standard methods using the Duhamel's principle and the fixed point argument, see for instance \cite{bachelot-2}, \cite{deldado}  and \cite{tartar}.  Now we suppose that for $(u_0,v_0)\in C_c^{\infty}$, the initial value problem (\ref{eq-dirac}-\ref{eq-dirac-initialv}) has a smooth solution $(u,v)\in C([0,T_*); H^1(R^1))\cap C^1([0,T_*); L^2(R^1))$ for some $T_*>0$. Then the support of $(u(x,t),v(x,t))$ is compact for each $t\in [0,T_*)$. To extend the solution across the time $t=T_*$ we only need to show that
 \begin{equation}\label{eq-H1-bound} \sup_{0\le t<T_*} \big(||u(t)||_{H^1} +||v(t)||_{H^1}\big)\le C^{\prime}(T_*)\exp(C^{\prime\prime}(T_*)) \end{equation}
 where  $C^{\prime}(\lambda)$ and $C^{\prime\prime}(\lambda)$ are polynomials of $\lambda$ independent of $T_*$.
To this end, we multiply the first equation of (\ref{eq-dirac}) by $\overline{u}$   and the second equation by $\overline{v}$, then
\begin{equation}\label{eq-dirac1}
\left\{\begin{array}{l}
(|u|^2)_t +(|u|^2)_x=mi(\overline{u}v-u\overline{v})+2\Re (i\overline{N_1}u), \\  (|v|^2)_t-(|v|^2)_x=-mi(\overline{u}v-u\overline{v})+2\Re (i\overline{N_2}v).
\end{array}\right.
\end{equation}
By Assumption (A2), we have
\[ (|u|^2+|v|^2)_t +(|u|^2-|v|^2)_x=0,\]
which leads to the following,
\begin{lemma}\label{lemma-conv-charge}
For any $t\in [0,T_*)$ there holds
\[ \int^{\infty}_{-\infty} (|u(t,x)|^2+|v(t,x)|^2)dx =\int^{\infty}_{-\infty} (|u_0(x)|^2+|v_0(x)|^2)dx.\]
\end{lemma}

Moreover, by Assumption (A1), we can find a constant $c>0$ such that for any $(u,v)$ there hold the following,
\begin{equation}\label{eq-nonlinearterm} |\overline{N_1}u|+ |\overline{N_2}v|\le c|u|^2|v|^2.\end{equation}
Then, (\ref{eq-dirac1}) implies that
\begin{equation}\label{eq-dirac10}
\left\{\begin{array}{l}
(|u|^2)_t +(|u|^2)_x\le r_0(t,x), \\  (|v|^2)_t-(|v|^2)_x\le r_0(t,x).
\end{array}\right.
\end{equation}
Here
\[ r_0(t,x)= m(|u|^2+|v|^2)(t,x)+ 2c|u(t,x)|^2|v(t,x)|^2.\]

As in \cite{bony} (see also \cite{ha-tzavaras} and \cite{zhou}), we define the followings  functional for the solution for any $t\in [0,T_*)$,
\[ Q_0(t)=\int\int_{x<y} |u(t,x)|^2|v(t,y)|^2 dxdy,\]
and \[ L_0(t)=\int^{\infty}_{-\infty} (|u(t,x)|^2+|v(t,x)|^2)dx,\]
\[ D_0(t)=\int^{\infty}_{-\infty} |u(t,x)|^2|v(t,x)|^2 dx.\]
$Q_0(t) $ is called a Bony functional \cite{ha-tzavaras} and is similar to the Glimm interaction potential.

\begin{lemma}\label{lemma-bony-functional} There exists  constants $\delta>0$ such that for the initial data satisfying $L_0(0)\le \delta$ there holds the following
\begin{eqnarray}\label{eq-bony-functional-ineq2}
\frac{d Q_0(t)}{dt }+D_0(t)\le 2m(L_0(0))^2
\end{eqnarray}
for $t\in [0,T_*)$.
Therefore,
\begin{eqnarray}\label{eq-bony-functional-ineq3}
 Q_0(t) + \int_0^t D_0(\tau) d\tau \le 2m(L_0(0))^2 t+ Q_0(0)\le 2m(L_0(0))^2 t+ (L_0(0))^2
\end{eqnarray}
for $t\in [0,T_*)$.
\end{lemma}
{\it Proof.} By (\ref{eq-dirac10}), we have
\begin{eqnarray*}
\frac{d Q_0(t)}{dt }  &\le& -\int\int_{x<y} (|u(t,x)|^2)_x |v(t,y)|^2 dxdy + \int\int_{x<y} |u(t,x)|^2 (|v(t,y)|^2)_y dxdy \\ & \, &+ \int\int_{x<y}\big(r_0(t,x)|v(t,y)|^2+ |u(t,x)|^2r_0(t,y) \big)  dxdy  \\
&\le& -2\int^{\infty}_{-\infty} |u(t,x)|^2|v(t,x)|^2 dx
\\
&\, & + \int^{\infty}_{-\infty} r_0(t,x)dx \int^{\infty}_{-\infty} |v|^2 dy + \int^{\infty}_{-\infty} |u|^2 dx\int^{\infty}_{-\infty} r_0(t,y)dy \\
 &\le & \big( -2+2L_0(t)c \big)\int^{\infty}_{-\infty} |u(x,t)|^2|v(x,t)|^2 dx + 2m(L_0(t))^2 \\
 &= & \big( -2+2L_0(0)c \big)D_0(t) + 2m(L_0(0))^2,
\end{eqnarray*}
where we also use Lemma \ref{lemma-conv-charge} in getting last equality.
Now, choose a constant $\delta>0$ such that
\[ -2+2\delta c<-1.\]Then, we can get the desired result. The proof is complete. $\Box$

Next we shall use the above estimates to get the $L^{\infty}$ bound of the solution for $0<t<T_*$. Denote \[\Omega(x_0,t_0)=\{(x,t): 0<t<t_0, \, x_0-(t_0-t)<x<x_0+(t_0-t)\},\] and
\[ \Gamma_R(x_0,x_0)=\{(x,t): 0<t<t_0, \, x=x_0+(t_0-t)\}, \]
\[ \Gamma_L(x_0,x_0)=\{(x,t): 0<t<t_0, \, x=x_0-(t_0-t)\}, \]
see Fig. \ref{fig-domain}.

\input dirac-fig.tex

\begin{lemma}\label{lemma-integral}
For any $t_0\in (0,T_*)$, there hold
\begin{eqnarray*}
\int_{\Gamma_R(x_0,t_0)} |u|^2 \le q(t_0), \\
\int_{\Gamma_L(x_0,t_0)} |v|^2 \le q(t_0),
\end{eqnarray*}
where $q(t_0)= t_0 (mL_0(0)+4c m(L_0(0))^2)+ 4c (L_0(0))^2+ L_0(0)$.
\end{lemma}
{\it Proof.} Integrating (\ref{eq-dirac1}) over $\Omega(x_0,t_0)$ yields
\begin{eqnarray*}
|\int\int_{\Omega(x_0,t_0)}\Big((|u|^2)_t +(|u|^2)_x\Big)dxdt| =\int\int_{\Omega(x_0,t_0)} 2|\Re \big(\overline{u}(m iv-iN_1)\big)| dxdt\\
\le \int\int_{\Omega(x_0,t_0)}\big( m(|u|^2+|v|^2)+2c |u|^2|v|^2\big) dxdt \\
\le mt_0 L_0(0) +  2c\int_0^{t_0}D(\tau) d\tau \\
\le  t_0 (mL_0(0)+4c m(L_0(0))^2)+ 4c (L_0(0))^2
\end{eqnarray*}
and
\begin{eqnarray*}
|\int\int_{\Omega(x_0,t_0)}\Big((|v|^2)_t -(|v|^2)_x\Big)dxdt| = 2|\Re \big(\overline{v}(m i u-iN_2)\big)| \\
\le mt_0 L_0(0) +  2c\int_0^{t_0}D(\tau) d\tau \\
\le  t_0 (mL_0(0)+4c m(L_0(0))^2)+ 4c (L_0(0))^2,
\end{eqnarray*}
where we use (\ref{eq-nonlinearterm}) and use Lemmas \ref{lemma-conv-charge} and \ref{lemma-bony-functional}.
Then the desired estimates follow by applying  Green formula to the above.
The proof is complete. $\Box$

\begin{lemma}\label{lemma-bound}
For $t\in [0,T_*)$, there holds the following,
\[ \sup_{ x\in R^1} \big( |u(x,t)|^2+|v(x,t)|^2 \big)\le \Big(\sup_{ x\in R^1} ( |u_0(x)|^2+|v_0(x)|^2 )+2mq(T_*)\Big)\exp\big(mT_*+2cq(T_*)\big) , \]
where $q(t)$ is given in Lemma \ref{lemma-integral}.
\end{lemma}
{\it Proof.}
 Using the characteristic method, by (\ref{eq-dirac10})  we have
\begin{eqnarray*}
\frac{d}{dt}|u(t,x_0+(t-t_0))|^2 \le \Big((m+2c|v|^2)|u|^2+m|v|^2\Big)\Big|_{(t,x_0+(t-t_0))}
\end{eqnarray*}
and
\begin{eqnarray*}
\frac{d}{dt}|v(t,x_0-(t-t_0))|^2\le \Big((m+2c|u|^2)|v|^2+m|u|^2\Big)\Big|_{(t,x_0-(t-t_0))},
\end{eqnarray*}
for any $x_0\in R^1$, $t_0\in [0,T_*)$.
Then,
\begin{eqnarray*}
|u(t_0,x_0)|^2 \le  ( |u_0(x_0-t_0)|^2 +m \int_0^{t_0} |v(\tau, x_0+(\tau-t_0))|^2d\tau) h_1(t_0),
 \\
|v(t_0,x_0)|^2 \le  (|v_0(x_0+t_0)|^2+m \int_0^{t_0}|u(\tau,x_0-(\tau-t_0))|^2d\tau) h_2(t_0),
\end{eqnarray*}
where
\begin{eqnarray*}
h_1(t_0)=\exp\big( mt_0+2c\int_0^{t_0}|v(x_0+(\tau-t_0),\tau)|^2 d\tau\big), \\
h_2(t_0)=\exp\big(mt_0+2c \int_0^{t_0}|u(x_0-(\tau-t_0),\tau)|^2 d\tau\big).
\end{eqnarray*}
 Therefore, these lead to the desired result by Lemma \ref{lemma-integral}. The proof is complete. $\Box$

 {\bf Proof of Theorem \ref{thm-1-existence}.} It suffices to get the $H^1$ bound of the solution on $[0,T_*)$. We differentiate Eqs. \ref{eq-dirac} with respect to $x$; then
 \begin{equation}\label{eq-dirac-H1}
\left\{ \begin{array}{l} i(u_{tx}+u_{xx})=-mv_x +\big(N_1(u,v)\big)_x, \\i( v_{tx}-v_{xx})=-mu_x+\big(N_2(u,v)\big)_x.
\end{array}
\right.
\end{equation}
For the nonlinear terms we have the following estimates by Lemma \ref{lemma-bound}:
\[ ||(N_1)_x||^2_{L^2(R^1)}+||(N_2)_x||^2_{L^2(R^1)}\le C_1q_1(T_*)(||u_x(t,\cdot)||^2_{L^2(R^1)}+||v_x(t,\cdot)||^2_{L^2(R^1)})  \]
for $t\in [0, T_*)$. Here  $C_1$ is a constant depending on $m$ and $\alpha_j$, $\beta_j$, $1\le j\le 5$; $q_1(T_*)=(1+q(T_*))^2\exp(2mT_*+2cq(T_*))$. This estimates, together with the energy method, we can get the desired bound as (\ref{eq-H1-bound}). Thus we can extend the solution across $t=T_*$. The proof is complete.$\Box$

\section{Convergence in $L^2$}

First we establish the estimates on the difference of smooth solutions. Let $(u^{\prime}, v^{\prime})$ be the global smooth solution to (\ref{eq-dirac}) taking $(u^{\prime}_0, v^{\prime}_0)\in C_c^{\infty}(R^1)$ as its initial data with $\int_{-\infty}^{\infty} (|u^{\prime}_0(x)|^2+|v^{\prime}_0(x)|^2)dx<\delta$.

Let $U=u-u^{\prime}$ and $V=v-v^{\prime}$. Then,
\begin{equation}\label{eq-dirac-difference1}
\left\{ \begin{array}{l} U_t+U_x=imV -i\Big(N_1(u,v)-N_1(u^{\prime},v^{\prime})\Big), \\V_t-V_x=imU-i\Big(N_2(u,v) -N_2(u^{\prime},v^{\prime})\Big),
\end{array}
\right.
\end{equation}
which leads to
\begin{equation}\label{eq-dirac-difference2}
\left\{ \begin{array}{l} (|U|^2)_t+(|U|^2)_x=\Re 2\{imV\overline{U} -i\Big(N_1(u,v)-N_1(u^{\prime},v^{\prime})\Big)\overline{U}\}, \\(|V|^2)_t-(|V|^2)_x=\Re 2\{imU\overline{V}-i\Big(N_2(u,v) -N_2(u^{\prime},v^{\prime})\Big)\overline{V}\}
\end{array}
\right.
\end{equation}
and
\begin{equation}
(|U|^2+|V|^2)_t+(|U|^2-|V|^2)_x=-R
\end{equation}
with
\[ R= \Re \{2i\Big(N_1(u,v)-N_1(u^{\prime},v^{\prime})\Big)\overline{U}+ 2i\Big(N_2(u,v) -N_2(u^{\prime},v^{\prime})\Big)\overline{V}\}.\]
For the nonlinear terms in righthand side in the above, we have the following.
\begin{lemma}\label{lemma-difference}
There exists a constant $c_*>0$ such that
\[ |\Re 2\{imV\overline{U} -i\Big(N_1(u,v)-N_1(u^{\prime},v^{\prime})\Big)\overline{U}\}| \le r_1(t,x),\]
\[ |\Re 2\{imU\overline{V}-i\Big(N_2(u,v) -N_2(u^{\prime},v^{\prime})\Big)\overline{V}\}|\le  r_1(t,x),\]
and
\[ |R|\le c_*r_2(t,x,x),\]
where $r_1(t,x)=m(|U(t,x)|^2+|V(t,x)|^2)+c_* r_2(t,x,x)$, and
\begin{eqnarray*}  r_2(t,x,y)= |U(t,x)|^2\Big( |v(t,y)|^2+|v^{\prime}(t,y)|^2\Big)  + \Big(|u(t,x)|^2+|u^{\prime}(t,x)|^2\Big)|V(t,y)|^2.
\end{eqnarray*}
\end{lemma}
{\it Proof.} Assumption (A1) implies that there exists a constant $c_*>0$ such that
\begin{eqnarray*}
|\Big(N_1(u,v)-N_1(u^{\prime},v^{\prime})\Big)\overline{U}|
&=&|N_1(u,v)-N_1(u^{\prime},v^{\prime})||U|  \\
& \le& \frac{c_*}{8}\big(|u-u^{\prime}||v|^2+|u^{\prime}||v-v^{\prime}||v|+|u^{\prime}||v^{\prime}||v-v^{\prime}| \big)|U| \\
& \le& \frac{c_*}{4} r_2(t,x,x),
\end{eqnarray*}
and
\[ |\Big(N_2(u,v)-N_2(u^{\prime},v^{\prime})\Big)\overline{V}| \le \frac{c_*}{4} r_2(t,x,x).\]
Thus, the above estimates give the desired results. The proof is complete.$\Box$

Now it follows from (\ref{eq-dirac-difference2})  that
\begin{equation}\label{eq-dirac-difference3}
\left\{ \begin{array}{l} (|U|^2)_t+(|U|^2)_x\le r_1(t,x), \\(|V|^2)_t-(|V|^2)_x\le r_1(t,x).
\end{array}
\right.
\end{equation}
And to deal with the nonlinear terms in (\ref{eq-dirac-difference1}), we define the following for any $t\ge 0$,
\[ L_1(t)=\int^{\infty}_{-\infty} (|U(t,x)|^2+|V(t,x)|^2)dx,\]
\[ Q_1(t)=\int\int_{x<y} r_2(t,x,y) dxdy,\]
\[ D_1(t)=\int^{\infty}_{-\infty} r_2(t,x,x) dx.\]

Let \[ r^{\prime}_0(t,x)=m(|u^{\prime}(t,x)|^2+|v^{\prime}(t,x)|^2)+c|u^{\prime}(t,x)|^2|v^{\prime}(t,x)|^2,\]
and
\[ Q_0^{\prime}(t)=\int\int_{x<y} |u^{\prime}(t,x)|^2|v^{\prime}(t,y)|^2 dxdy,\]
and \[ L_0(t)=\int^{\infty}_{-\infty} (|u^{\prime}(t,x)|^2+|v^{\prime}(t,x)|^2)dx,\]
\[ D_0(t)=\int^{\infty}_{-\infty} |u^{\prime}(t,x)|^2|v^{\prime}(t,x)|^2 dx.\]

\begin{lemma}\label{lemma-functinal-difference}
There exist constants $\delta>0$ and $K>0$ such that if $L_0(0)<\delta$ and $L_0^{\prime}(0)<\delta$ then
\begin{equation}\label{eq-bony-functional-4} \frac{d}{dt}(L_1(t)+KQ_1(t))+D_1(t) \le \Big(2mL_0(0)+2mL_0^{\prime}(0)+cD_0(t)+cD^{\prime}_0(t)\Big)L_1(t)
\end{equation}
for $t\ge 0$. Therefore
\begin{equation}\label{eq-bony-functional-4-1}
 L_1(t)+KQ_1(t)\le (L_1(0)+KQ_1(0)) \exp(h_3(t))\end{equation} and
 \begin{equation}\label{eq-bony-functional-5}
\int_0^t D_1(\tau)d\tau\le (L_1(0)+KQ_1(0)) (1+(4m\delta+\delta^2+2m\delta^2 t)\int^t_0\exp(h_3(\tau))d\tau)
\end{equation}
for for $t\ge 0$, where
\[h_3(t)=2mL_0(0)t+2mL_0^{\prime}(0)t+\int_0^t(cD_0(\tau)+cD^{\prime}_0(\tau))d\tau.\]
\end{lemma}
{\it Proof.} (\ref{eq-dirac-difference3}) yields that
%\begin{equation}\label{eq-bony-functional-difference1}
\[ \frac{d}{dt}L_1(t) \le 2c_* D_1(t),\]
%\end{equation}
and
\begin{eqnarray*}
\frac{d}{dt}Q_1(t) &\le& -2D_1(t) + \int\int_{x<y} r_1(t,x) (|v(t,y)|^2+|v^{\prime}(t,y)|^2)dxdy \\
 &\,& + \int\int_{x<y} |U(t,x)|^2(r_0(t,y)+r_0^{\prime}(t,y))dxdy \\
 &\,& +\int\int_{x<y} (|u(t,x)|^2+|u^{\prime}(t,x)|^2)r_1(t,x)dxdy  \\
 &\,& +\int\int_{x<y} (r_0(t,x)+r_0^{\prime}(t,x))|V(t,y)|^2 dxdy  \\
 &\le& -2D_1(t)+ (mL_1(t)+c_* D_1(t))(L_0(t)+L_0^{\prime}(t)) \\
 &\,& + (mL_0(t)+mL_0^{\prime}(t)+cD_0(t)+cD^{\prime}_0(t))L_1(t) \\
 &\le& \Big(-2+c_*(L_0(0)+L_0^{\prime}(0))\Big)D_1(t) \\ &\,&+ \Big(2mL_0(0)+2mL_0^{\prime}(0)+cD_0(t)+cD^{\prime}_0(t)\Big)L_1(t),
\end{eqnarray*}
where Lemma \ref{lemma-conv-charge} is used in the last inequality. These estimates give (\ref{eq-bony-functional-4})  for small $L_0(0)<\delta$ and $L_0^{\prime}(0)<\delta$ and for Large $K$, where $\delta=\frac{1}{4c_*}$. (\ref{eq-bony-functional-4-1}) is a consequence of (\ref{eq-bony-functional-4}).

To prove (\ref{eq-bony-functional-5}), we integrate (\ref{eq-bony-functional-4}) over $[0,t]$, then
\begin{eqnarray*}
\int_0^t D_1(\tau)d\tau &\le& L_1(0)+KQ_1(0)+\int_0^t \Big(4m\delta+cD_0(\tau)+cD^{\prime}_0(\tau)\Big)L_1(\tau)d\tau \\
 &\le& (L_1(0)+KQ_1(0))+(4m\delta+\delta^2+2m\delta^2 t)\int_0^t L_1(\tau)d\tau \\
 &\le& (L_1(0)+KQ_1(0)) (1+(4m\delta+\delta^2+2m\delta^2 t)\int^t_0\exp(h_3(\tau))d\tau),
\end{eqnarray*}
where Lemma \ref{lemma-bony-functional} and (\ref{eq-bony-functional-4-1}) are used.
The proof is complete. $\Box$

\begin{lemma}\label{lemma-h3(t)}
Let $h_3(t)$ be given by Lemma \ref{lemma-functinal-difference}. Then for $t\ge 0$,
\begin{eqnarray*}
 h_3(t)\le 2m\Big(L_0(0)+L_0^{\prime}(0)+(L_0(0))^2+ (L_0^{\prime}(0))^2 \Big)t+ (L_0(0))^2+(L_0^{\prime}(0))^2.
\end{eqnarray*}
\end{lemma}

\begin{proposition}\label{prop-difference} Let $(u,v)$ and $(u^{\prime},v^{\prime})$ be two classical solutions to \ref{eq-dirac} with initial data $(u_0,v_0)$ and $(u_0^{\prime}, v^{\prime})$ respectively, and suppose that $||u_0||_{L^2(R^1)}^2+||v_0||_{L^2(R^1)}^2<\delta $ and $||u_0^{\prime}||_{L^2(R^1)}^2+||v_0^{\prime}||_{L^2(R^1)}^2<\delta $. Then, there exist constants, $c_1$, $c_2$ and $c_3$, depending on the $\delta$, such that for $t\ge 0$,
\[ ||u(t)-u^{\prime}(t)||_{L^2(R^1)}^2+||v(t)-v^{\prime}(t)||_{L^2(R^1)}^2\le (||u_0-u^{\prime}_0||_{L^2(R^1)}^2 +||v_0-v^{\prime}_0||_{L^2(R^1)}^2)h_4(t),\]
where $h_4(t)=c_1\exp(c_2t+c_3)$.
\end{proposition}
{\it Proof.} Indeed, we have
\[ Q_1(0)\le \int_{-\infty}^{\infty}( |U(0,x)|^2+|V(0,x)|^2)dx (L_0(0)+L_0^{\prime}(0)),\]
which, together with (\ref{eq-bony-functional-4-1}), gives the proof of this proposition. Thus the proof is complete.$\Box$

{\bf Proof of Theorem \ref{thm-2-convergence}} By Proposition \ref{prop-difference}, we have
\begin{eqnarray*}
& ||u_k(t)-u_j(t)||_{L^2(R^1)}^2+||v_k(t)-v_j(t)||_{L^2(R^1)}^2 \\ &\le (||u_{k0}-u_{j0}||_{L^2(R^1)}^2 +||v_{k0}-v_{j0}||_{L^2(R^1)}^2)h_4(t),
\end{eqnarray*}
which implies that there exists a unique $(u_{\infty},v_{\infty})\in C([0,+\infty),L^(R^2))$ such that
\[ ||u_k-u_{\infty}||_{C([0,T],L^(R^2))}+||v_k-v_{\infty}||_{C([0,T],L^(R^2))}\to 0\]
for any $T\ge 0$.

Now to prove that $(u_{\infty},v_{\infty})$ is a weak solution of (\ref{eq-dirac}), let
\[ U_{k,j}=u_k-u_j,\quad V_{k,j}=v_k-v_j,\]
and
\[ r_{k,j}= |U_{k,j}(t,x)|^2\Big( |v_k(t,x)|^2+|v_j(t,x)|^2\Big)  + \Big(|u_k(t,x)|^2+|u_j(t,x)|^2\Big)|V_{k,j}(t,x)|^2.
\]
By (\ref{eq-bony-functional-5}),  for any $T>0$, there exists a constant $C_1(T)>0$ such that for $l=1,2$
\begin{eqnarray*}
\int_{-\infty}^{\infty}|N_l(u_k,v_k)-N_l(u_j,v_j)| dx &\le& c_*\int_{-\infty}^{\infty} r_{k,j} dx \\
&\le & C_1(T)c_* (||U_{k,j}(t)||_{L^2(R^1)}^2+||V_{k,j}(t)||_{L^2(R^1)}^2) \\ 
&\le & C_1(T)c_* (||U_{k,j}(0)||_{L^2(R^1)}^2+||V_{k,j}(0)||_{L^2(R^1)}^2),
\end{eqnarray*}
which leads to the strong convergence of the nonlinear terms. Therefore, $(u_{\infty},v_{\infty})$ is a weak solution of (\ref{eq-dirac}). The proof of Theorem \ref{thm-2-convergence} is complete. $\Box$

\section*{ Acknowledgement}

This work was partially  supported by NSFC Project 11031001 and 11121101, by the 111 Project
B08018 and by Ministry of Education of China.

\end{document}

%% file: dirac-fig.tex
\begin{figure}[h]
\begin{center}
\unitlength=1mm
%\thicklines
\begin{picture}(130,75)
\put(40,40){$\Omega(x_0,t_0)$}
\put(10,43){$\Gamma_L(x_0,t_0)$}
\put(60,45){$\Gamma_R(x_0,t_0)$}
\put(85,20){\line(-1,1){40}}
\put(5,20){\line(1,1){40}} \put(43,63){$(x_0,t_0)$}
\put(0,20){\vector(1,0){100}}\put(102,20){$x$}
\put(0,16){$(x_0-t_0, 0)$ } \put(78,16){$(x_0+t_0,0)$}
\end{picture}
\caption{Domain $\Omega(x_0,t_0)$}\label{fig-domain}
\end{center}
\end{figure}
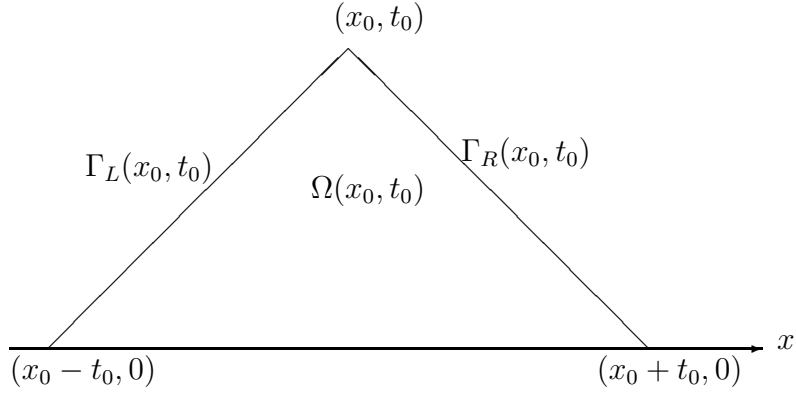

%% file: dirac-massive20130406.bbl
\begin{thebibliography}{9}



\bibitem{bachelot} A. Bachelot, {\it Global existence of large amplitude solutions for nonlinear massless Dirac equation,} Portucaliae Math.{\bf 46} Fasc. Supl. (1989), 455-473.

\bibitem{bachelot-2} A. Bachelot, {\it Global Cauchy problem for semilinear hyperboloc systems with nonlocal interactions. Applications to Dirac equations,} J. Math. Pures. Appl. {\bf 86} (2006), 201-236.

\bibitem{bony} J. M. Bony, {\it Solution globales born\'{e}es pour les mod\`{e}les discrets de l'\'{e}quation de Boltzmann, en dimension 1 d'espace,} Acte des Journ\'{e}es E.D.P. \`{a} Saint-Jean-de-Monts, Ecole Polytechnique (1987), XVI-1-XVI-9.

\bibitem{bournaveas-zouraris} N. Bournaveas and G. E. Zouraris, {\it Theory and numerical approximations for a nonlinear $1+1$ Dirac system,}  ESAIM: Math. Model. Num. Analysis {\bf 46} (4) (2012), 841-874.

\bibitem{bressan} A. Bressan, {\it Hyperbolic systems of conservation laws:
 The one-dimensional Cauchy problem,} Oxford University Press
 Inc., New York, 2000.

\bibitem{cacciafesta1} F. Cacciafesta, {\it Global small solutions to the critical radial Dirac equation with potential,} Nonlinear Analysis {\bf 74} (2011), 6060-6073.

\bibitem{cacciafesta2} F. Cacciafesta, {\it  The cubic nonlinear Dirac equation,} arXiv:1301.6867v1 [math.AP] 29 Jan 2013.

\bibitem{candy} T. Candy, {\it Global existence for an $L^2$ critical nonlinear Dirac equation in one dimension,}  Adv. Differential Equations {\bf 16} No. 7-8 (2011), 643-666.


\bibitem{dafermos} C. M. Dafermos, {\it Hyperbolic Conservation Laws in Continuum Physics}, Springer-Verlag, Berlin, 2010.

\bibitem{deldado} V. Delgado, {\it Global solution of the Cauchy problem for the (classical) coupled Maxwell-Dirac and other nonlinear Dirac equations,} Proc. Amer. Math. Soc. {\bf 69}(2)(1978), 289-296.

\bibitem{dias-figueira1} J. P. Dias and M. Figueira, {\it Sur l'existence d'une solution globale pour une equation de Dirac non lin\'{e}aire avec masse nulle,} C.R. Acad. Sci. Paris, S\'{e}rie I, {\bf 305} (1987), 469-472.

\bibitem{dias-figueira2} J. P. Dias and M. Figueira, {\it Remarque sur le probl\`{e}me de Cauchy pour une equation de Dirac non lin\'{e}aire avec masse nulle,} Portucaliae Math. {\bf 45}(4) (1988), 327-335.

\bibitem{escobedo} M. Escobedo and L. Vega, {\it A semilinear Dirac equation in $H^s(R^3)$ for $s>1$,} SIAM J. Math. Anal. {\bf 28}(2) (1997), 338-362.

\bibitem{Esteban-Lewin-Sere} M. J. Esteban, M. Lewin and E. S\'{e}r\'{e}, {\it Variational methods in relativistic quantum mechanics,} Bulletin A.M.S. {\bf 45} No. 4 (2008), 535-598.

\bibitem{ha-tzavaras} S-Y Ha and A. E. Tzavaras, {\it Lyapunov functionals and $L^1$-stability for discrete velocity Boltzmann equations,} Commun. Math. Phys. {\bf 239} (2003), 65-92.

\bibitem{glimm} J. Glimm, {\it Solution in the large
for nonlinear systems of conservation laws,} Comm.
Pure Appl. Math. {\bf 18}(1965), 695-715.

\bibitem{gross-neveu} D. J. Gross and A. Neveu, {\it Dynamical symmetry breaking in asymptotically free field theories,} Phys. Rev. D 10 (1974), 3235-3253.

\bibitem{huh} H. Huh, {\it Global strong solution to the Thirring model in critical space,} J. Math. Anal. Appl. {\bf 381} (2011), 513-520.

\bibitem{machihara} S. Machihara and T. Omoso, {\it The explicit solutions to the nonlinear Dirac equation and Dirac-Klein-Gordon equation,} Ricerche Mat. {\bf 56} (1)(2007), 19-30.

\bibitem{pelinovsky} D. Pelinovsky, {\it Survey on global existence in the nonlinear Dirac equations in one dimension,} arXiv:1011.5925vl [math-ph] 26 Nov 2010.

\bibitem{selberg} S. Selberg and A. Tesfahun, {\it Low regularity well-posedness for some nonlinear Dirac equations in one space dimension,} Differential Integral Equations {\bf 23}(3-4)(2010), 265-278.

\bibitem{tartar} L. Tartar, {\it From Hyperbolic Systems to Kinetic Theory,} Springer 2008.

\bibitem{thirring} W.E. Thirring, {\it A soluble relativistic field theory,} Ann. Phys. {\bf 3}(1958), 91-112.

\bibitem{zhang} Y. Zhang, {\it Global strong solution to a nonlinear Dirac type equation in one dimension,} Nonlinear Analysis: Theory, Method and Applications {\bf 80} (2013), 150-155.

\bibitem{zhou} Y. Zhou, {\it Uniqueness of weak solutions in 1+1 dimensional wave maps,} Math. Z. {\bf 232}(1999), 707-719.


\end{thebibliography}
